\documentclass[12pt]{amsart}
\usepackage{amssymb}

\newcommand{\Bgp}{{\Z^\N}}

\long\def\forget#1\forgotten{}
\newcommand{\issuenumber}{29}
\newcommand{\issuemonth}{December}
\newcommand{\issueyear}{2009}

\newcommand{\ed}{
\newpage

\section{Unsolved problems from earlier issues}

\begin{issue}
Is $\binom{\Omega}{\Gamma}=\binom{\Omega}{\Tau}$?
\end{issue}

\begin{issue}
Is $\ufin(\cO,\Omega)=\sfin(\Gamma,\Omega)$?
And if not, does $\ufin(\cO,\Gamma)$ imply
$\sfin(\Gamma,\Omega)$?
\end{issue}

\stepcounter{issue}

\begin{issue}
Does $\sone(\Omega,\Tau)$ imply $\ufin(\Gamma,\Gamma)$?
\end{issue}

\begin{issue}
Is $\fp=\fp^*$? (See the definition of $\fp^*$ in that issue.)
\end{issue}

\begin{issue}
Does there exist (in ZFC) an uncountable set satisfying $\sfin(\cB,\cB)$?
\end{issue}

\stepcounter{issue}

\begin{issue}
Does $X \nin \NON(\cM)$ and $Y\nin\mathsf{D}$ imply that
$X\cup Y\nin \COF(\cM)$?
\end{issue}

\begin{issue}[CH]
Is $\split(\Lambda,\Lambda)$ preserved under finite unions?
\end{issue}

\begin{issue}
Is $\cov(\cM)=\fo$? (See the definition of $\fo$ in that issue.)
\end{issue}

\begin{issue}\label{b}
Does $\sone(\Gamma,\Gamma)$ always contain an element of cardinality $\fb$?
\end{issue}

Problem \ref{b} is solved. See Section \ref{pcf} above.

\begin{issue}
Could there be a Baire metric space $M$ of weight $\aleph_1$ and a partition
$\mathcal{U}$ of $M$ into $\aleph_1$ meager sets where for each ${\mathcal U}'\subset\mathcal U$,
$\bigcup {\mathcal U}'$ has the Baire property in $M$?
\end{issue}

\stepcounter{issue} 

\begin{issue}
Does there exist (in ZFC) a set of reals $X$ of cardinality $\fd$ such that all
finite powers of $X$ have Menger's property $\sfin(\cO,\cO)$?
\end{issue}

\begin{issue}
Can a Borel non-$\sigma$-compact group be generated by a Hurewicz subspace?
\end{issue}

\begin{issue}[MA]
Is there $X\sbst\bbR$ of cardinality continuum, satisfying $\sone(\BO,\BG)$?
\end{issue}

\begin{issue}[CH]
Is there a totally imperfect $X$ satisfying $\ufin(\cO,\Gamma)$
that can be mapped continuously onto $\Cantor$?
\end{issue}

\begin{issue}[CH]
Is there a Hurewicz $X$ such that $X^2$ is Menger but not Hurewicz?
\end{issue}

\begin{issue}
Does the Pytkeev property of $C_p(X)$ imply that $X$ has Menger's property?
\end{issue}

\begin{issue}
Does every hereditarily Hurewicz space satisfy $\sone(\BG,\BG)$?
\end{issue}

\begin{issue}[CH]
Is there a Rothberger-bounded $G\le\Bgp$ such that $G^2$ is not Menger-bounded?
\end{issue}

\begin{issue}
Let $\cW$ be the van der Waerden ideal.
Are $\cW$-ultrafilters closed under products?
\end{issue}

\begin{issue}
Is the $\delta$-property equivalent to the $\gamma$-property $\binom{\Omega}{\Gamma}$?
\end{issue}

\stepcounter{issue}

\stepcounter{issue}

\general\end{document}}

\newcommand{\Cantor}{{\{0,1\}^\N}}

\newcommand{\fb}{\mathfrak{b}}

\newcommand{\fc}{\mathfrak{c}}
\newcommand{\fd}{\mathfrak{d}}
\newcommand{\fp}{\mathfrak{p}}

\newcommand{\NON}{{\mathsf   {NON}}}
\newcommand{\COF}{{\mathsf   {COF}}}

\newcommand{\E}{\exists}

\newcommand{\cA}{\mathcal{A}}
\newcommand{\cM}{\mathcal{M}}

\newcommand{\op}{\operatorname}
\newcommand{\cov}{\mathsf{cov}}

\newcommand{\cof}{\mathsf{cof}}

\newcommand{\bbR}{\mathbb{R}}

\newcommand{\fo}{\mathfrak{od}}

\newcommand{\w}{\omega}

\renewcommand{\split}{\mathsf{Split}}
\newcommand{\bq}{\begin{quote}}
\newcommand{\eq}{\end{quote}}
\newcommand{\cO}{\mathcal{O}}
\newcommand{\cB}{\mathcal{B}}
\newcommand{\BG}{\cB_\Gamma}

\newcommand{\BO}{\cB_\Omega}

\newcommand{\sone}{\mathsf{S}_1}    \newcommand{\sfin}{\mathsf{S}_\mathrm{fin}}

\newcommand{\ufin}{\mathsf{U}_\mathrm{fin}}

\newcommand{\nin}{\not\in}

\newcommand{\cU}{\mathcal{U}}

\newcommand{\cW}{\mathcal{W}}

\newcommand{\NN}{{\N^\N}}
\newcommand{\N}{\mathbb{N}}
\newcommand{\bbN}{\mathbb{N}}
\newcommand{\Z}{\mathbb{Z}}

\newcommand{\sbst}{\subseteq}
\newcommand{\by}[2]{\par\hfill\emph{#1}, #2}
\newcommand{\nby}[1]{\par\hfill\emph{#1}}
\newcommand{\Tau}{\mathrm{T}}
\newcommand{\CE}{\textsc{CE}}

\newtheorem{thm}{Theorem}[section]
\newcommand{\bthm}{\begin{thm}} \newcommand{\ethm}{\end{thm}}
\newtheorem{prop}[thm]{Proposition}
\newcommand{\bprp}{\begin{prop}} \newcommand{\eprp}{\end{prop}}
\newtheorem{fact}[thm]{Fact}
\newcommand{\bfct}{\begin{fact}} \newcommand{\efct}{\end{fact}}
\newtheorem{prob}[thm]{Problem}
\newcommand{\bprb}{\begin{prob}} \newcommand{\eprb}{\end{prob}}
\newtheorem{lem}[thm]{Lemma}
\newcommand{\blem}{\begin{lem}} \newcommand{\elem}{\end{lem}}
\newtheorem{claim}[thm]{Claim}
\newcommand{\bclm}{\begin{claim}} \newcommand{\eclm}{\end{claim}}
\newtheorem{cor}[thm]{Corollary}
\newcommand{\bcor}{\begin{cor}} \newcommand{\ecor}{\end{cor}}
\newtheorem{conj}[thm]{Conjecture}
\newcommand{\bcnj}{\begin{conj}} \newcommand{\ecnj}{\end{conj}}
\theoremstyle{definition}
\newtheorem{defn}[thm]{Definition}
\newcommand{\bdfn}{\begin{defn}} \newcommand{\edfn}{\end{defn}}
\theoremstyle{remark}
\newtheorem{rem}[thm]{Remark}
\newcommand{\brem}{\begin{rem}} \newcommand{\erem}{\end{rem}}
\newtheorem{cnv}[thm]{Convention}
\newcommand{\bcnv}{\begin{cnv}} \newcommand{\ecnv}{\end{cnv}}
\newtheorem{exam}[thm]{Example}
\newcommand{\bexm}{\begin{exam}} \newcommand{\eexm}{\end{exam}}
\newtheorem{issue}{Issue}

\newcommand{\bpf}{\begin{proof}} \newcommand{\epf}{\end{proof}}
\newcommand{\be}{\begin{enumerate}}
\newcommand{\ee}{\end{enumerate}}
\newcommand{\bi}{\begin{itemize}}
\newcommand{\ei}{\end{itemize}}
\newcommand{\itm}{\item}

\newcommand{\general}{\small\vfill\par\noindent\hrulefill\par
\noindent\textbf{Previous issues.} The previous issues of this
bulletin are available online at\\
\texttt{http://front.math.ucdavis.edu/search?\&t=\%22SPM+Bulletin\%22}
\\[0.1cm]
\textbf{Contributions.} Announcements, discussions, and open problems should be emailed
to \texttt{tsaban@math.biu.ac.il}\\[0.1cm]
\textbf{Subscription.}
To receive this bulletin (free) to your e-mailbox, e-mail us.
}

\newcommand{\arXiv}[5]{\subsection{#2}{#4}\par\hfill{\arx{#1}}\par\hfill\emph{#3}}

\newcommand{\arx}[1]{\texttt{http://arxiv.org/abs/#1}}
\newcommand{\url}[1]{\bq\texttt{#1}\eq}
\newcommand{\online}[1]{The paper is available online at \url{#1}}

\title[$\mathcal{SPM}$ Bulletin \textbf{\issuenumber} (\issuemonth{} \issueyear)]{%
$\mathcal{SPM}$ Bulletin\\[0.5cm]
Issue number \issuenumber: \issuemonth{} \issueyear{} \CE{}}

\begin{document}
\maketitle

\tableofcontents

\section{Editor's note}

The \emph{Problem of the Issue} posed in Issue 11 (see at the end of the present Issue)
was solved in Announcement \ref{pcf} below.

The Borel version of the problem from Issue 23 (see below)
was posed by Miller in his plenary lecture at the \emph{IIIrd Workshop on Coverings, Selections, and Games in Topology}.
This latter problem was also solved recently, by my student Tal Orenshtein and I.
We hope to have a paper describing this solution ready during 2010.

Several additional problems from the \emph{Open Problems in Topology} book chapter on SPM
are being solved in several works in progress by several authors.
This makes the civilian year 2009 very satisfactory in this respect.

Have a peaceful and fruitful 2010.

\medskip

\by{Boaz Tsaban}{tsaban@math.biu.ac.il}

\hfill \texttt{http://www.cs.biu.ac.il/\~{}tsaban}

\section{Research announcements in ``core'' SPM}

\noindent\emph{Note.} The division between the present section and the next one
is somewhat artificial, but perhaps still useful.

\arXiv{0909.3663}
{Measurable cardinals and the cardinality of Lindel\"of spaces}
{Marion Schepers}
{If it is consistent that there is a measurable cardinal, then it is
consistent that all points $G_\delta$ Rothberger spaces have ``small'' cardinality.}

\arXiv{0909.5004}
{Topological games and covering dimension}
{Liljana Babinkostova}
{We consider a natural way of extending the Lebesgue covering dimension to
various classes of infinite dimensional separable metric spaces.}

\arXiv{0909.5645}
{Menger's and Hurewicz's Problems: Solutions from ``The Book'' and refinements}
{Boaz Tsaban}
{We provide simplified solutions of Menger's and Hure\-wicz's problems
and conjectures, concerning generalizations of $\sigma$-compactness.
The reader who is new to this field will find a self-contained treatment
in Sections 1, 2, and 5.
\par
Sections 3 and 4 contain new results, based on the mentioned simplified solutions.
The main new result is that there is a set of reals $X$ of cardinality $\fb$, which
has the following property:
\begin{quote}
Given point-cofinite covers $\cU_1,\cU_2,\dots$ of $X$, there are for each $n$ sets $U_n,V_n\in\cU_n$,
such that each member of $X$ is contained in all but finitely many of the sets
$U_1\cup V_1,U_2\cup V_2,\dots$
\end{quote}
This property is strictly stronger than Hurewicz's covering property, and
by a result of Miller and the present author,
one cannot prove the same result if we are only allowed to pick one set from each $\cU_n$.}

\arXiv{0910.4063}
{Point-cofinite covers in the Laver model}
{Arnold W. Miller and Boaz Tsaban}
{\label{pcf}Let $\sone(\Gamma,\Gamma)$ be the statement: For each sequence of point-cofinite open covers
(i.e., $\gamma$-covers),
one can pick one element from each cover and obtain a point-cofinite cover.
$\fb$ is the minimal cardinality of a set of reals not satisfying $\sone(\Gamma,\Gamma)$.
We prove the following assertions:
\be
\itm If there is an unbounded tower, then there are sets of reals of cardinality $\fb$, satisfying $\sone(\Gamma,\Gamma)$.
\itm It is consistent that all sets of reals satisfying $\sone(\Gamma,\Gamma)$ have cardinality smaller than $\fb$.
\ee
These results can also be formulated as dealing with Arhangel'ski\u{\i}'s property $\alpha_2$
for spaces of continuous real-valued functions.
\par
The main technical result is that in Laver's model,
each set of reals of cardinality $\fb$ has an unbounded
Borel image in the Baire space $\NN$.}

\subsection{Projective versions of selection principles}

All spaces are assumed to be Tychonoff. A space $X$ is called projectively $P$
(where $P$ is a topological property) if every continuous second countable image of $X$ is $P$.
Characterizations of projectively Menger spaces $X$ in terms of continuous mappings $f : X\to\bbR^\w$,
of Menger base property with respect to separable pseudometrics and a selection principle
restricted to countable covers by cozero sets are given. If all finite powers of X
are projectively Menger, then all countable subspaces of $C_p(X)$ have countable fan tightness.
The class of projectively Menger spaces contains all Menger spaces as well as all
$\sigma$-pseudocompact spaces, and all spaces of cardinality less than $\fd$. Projective versions
of Hurewicz, Rothberger and other selection principles satisfy properties similar to the
properties of projectively Menger spaces, as well as some specific properties. Thus, $X$ is
projectively Hurewicz if and only if $C_p(X)$ has the Monotonic Sequence Selection Property in the
sense of Scheepers; $\beta X$ is Rothberger iff $X$ is pseudocompact and projectively Rothberger.
Embeddability of the countable fan space $V_\w$ into $C_p(X)$ or $C_p(X,2)$ is characterized in
terms of projective properties of X.

To appear in \emph{Topology and its Applications}:
\url{http://dx.doi.org/10.1016/j.topol.2009.12.004}

\nby{Maddalena Bonanzinga, Filippo Cammaroto, Mikhail Matveev}

\section{Additional research announcements}\label{RA}

\arXiv{0908.0475}
{Some Ramsey theorems for finite $n$-colorable and $n$-chromatic graphs}
{L. Nguyen Van Th\'e}
{Given a fixed integer $n$, we prove Ramsey-type theorems for the classes of
all finite ordered $n$-colorable graphs, finite $n$-colorable graphs, finite
ordered $n$-chromatic graphs, and finite $n$-chromatic graphs.}

\arXiv{0908.1216}
{Uniform convexity and the splitting problem for selections}
{Maxim V. Balashov and Du\v{s}an Repov\v{s}}
{We continue to investigate cases when the Repov\v{s}-Semenov splitting
problem for selections has an affirmative solution for continuous set-valued
mappings. We consider the situation in infinite-dimensional uniformly convex
Banach spaces. We use the notion of Polyak of uniform convexity and modulus of
uniform convexity for arbitrary convex sets (not necessary balls). We study
general geometric properties of uniformly convex sets. We also obtain an
affirmative solution of the splitting problem for selections of certain
set-valued mappings with uniformly convex images.}

\arXiv{0908.1544}
{The Solecki dichotomy for functions with analytic graphs}
{Janusz Pawlikowski, Marcin Sabok}
{A dichotomy discovered by Solecki says that a Baire class 1 function from a
Souslin space into a Polish space either can be decomposed into countably many
continuous functions, or else contains one particular function which cannot be
so decomposed. In this paper we generalize this dichotomy to arbitrary
functions with analytic graphs. We provide a ``classical'' proof, which uses only
elementary combinatorics and topology.}

\arXiv{0908.1605}
{A co-analytic maximal set of orthogonal measures}
{Vera Fischer, Asger Tornquist}
{We prove that if $V=L$ then there is a $\Pi^1_1$ maximal orthogonal (i.e.
mutually singular) set of measures on Cantor space. This provides a natural
counterpoint to the well-known Theorem of Preiss and Rataj that no analytic set
of measures can be maximal orthogonal.}

\arXiv{0908.2790}
{Model theory of operator algebras I: Stability}
{Ilijas Farah, Bradd Hart, David Sherman}
{Several authors have considered whether the ultrapower and the relative
commutant of a $C^*$-algebra or II$_1$ factor depend on the choice of the
ultrafilter. We show that the negative answer to each of these questions is
equivalent to the Continuum Hypothesis, extending results of Ge-Hadwin and the
first author.}

\arXiv{0908.2225}
{Characterizing meager paratopological groups}
{T. Banakh, I. Guran and A. Ravsky}
{We prove that a Hausdorff paratopological group $G$ is meager if and only if
there are a nowhere dense subset $A$ of $G$ and a countable subset $C$ in $G$ such that
$CA=G=AC$.}

\arXiv{0908.2228}
{The topological structure of direct limits in the category of uniform spaces}
{Taras Banakh}
{Let $(X_n)_{n}$ be a sequence of uniform spaces such that each space $X_n$ is
a closed subspace in $X_{n+1}$. We give an explicit description of the topology
and uniformity of the direct limit $u-lim X_n$ of the sequence $(X_n)$ in the
category of uniform spaces. This description implies that a function $f:u-lim
X_n\to Y$ to a uniform space $Y$ is continuous if for every $n$ the restriction
$f|X_n$ is continuous and regular at the subset $X_{n-1}$ in the sense that for
any entourages $U\in\cU_Y$ and $V\in\cU_X$ there is an entourage $V\in\cU_X$ such
that for each point $x\in B(X_{n-1},V)$ there is a point $x'\in X_{n-1}$ with
$(x,x')\in V$ and $(f(x),f(x'))\in U$. Also we shall compare topologies of
direct limits in various categories.}

\arXiv{0908.1942}
{Orthonormal bases of Hilbert spaces}
{Ilijas Farah}
{I prove that a Hilbert space has the property that each of its dense (not
necessarily closed) subspaces contains an orthoormal basis if and only if it is
separable.}

\arXiv{0908.1943}
{A dichotomy for the Mackey Borel structure}
{Ilijas Farah}
{We prove that the equivalence of pure states of a separable C*-algebra is
either smooth or it continuously reduces $[0,1]^{\bbN}/\ell_2$ and it therefore
cannot be classified by countable structures. The latter was independently
proved by Kerr--Li--Pichot by using different methods.
We also give some remarks on a 1967 problem of Dixmier.}

\arXiv{0909.4563}
{Minimal Size of Basic Families}
{Ziqin Feng and Paul Gartside}
{A family $\cB$ of continuous real-valued functions on a space $X$ is said
to be {\sl basic} if every $f \in C(X)$ can be represented $f = \sum_{i=1}^n
g_i \circ \phi_i$ for some $\phi_i \in \cB$ and $g_i \in C(\bbR)$ ($i=1, ...,
n$). Define $\op{basic} (X) = \min \{|\cB| : \cB$ is a basic family for $X\}$.
If $X$ is separable metrizable $X$ then either $X$ is locally compact and
finite dimensional, and $\op{basic} (X) < \aleph_0$, or $\op{basic} (X) =
\mathfrak{c}$. If $K$ is compact and either $w(K)$ (the minimal size of a basis
for $K$) has uncountable cofinality or $K$ has a discrete subset $D$ with
$|D|=w(K)$ then either $K$ is finite dimensional, and $\op{basic} (K) = \cof
([w(K)]^{\aleph_0}, \subseteq)$, or $\op{basic} (K) = |C(K)|=w(K)^{\aleph_0}$.}

\arXiv{0909.4561}
{On Hilbert's 13th Problem}
{Ziqin Feng and Paul Gartside}
{Every continuous function of two or more real variables can be written as the
superposition of continuous functions of one real variable along with addition.}

\arXiv{0909.5668}
{Metastability and the Furstenberg-Zimmer Tower II: Polynomial and
 Multidimensional Szemeredi's Theorem}
{Henry Towsner}
{The Furstenberg-Zimmer structure theorem for $\mathbb{Z}^d$ actions says that
every measure-preserving system can be decomposed into a tower of primitive
extensions. Furstenberg and Katznelson used this analysis to prove the
multidimensional Szemer\'edi's theorem, and Bergelson and Liebman further
generalized to a polynomial Szemer\'edi's theorem. Beleznay and Foreman showed
that, in general, this tower can have any countable height. Here we show that
these proofs do not require the full height of this tower; we define a weaker
combinatorial property which is sufficient for these proofs, and show that it
always holds at fairly low levels in the transfinite construction
(specifically, $\omega^{\omega^{\omega^\omega}}$).}

\arXiv{0910.0279}
{Cofinitary Groups and Other Almost Disjoint Families of Reals}
{Bart Kastermans}
{We study two different types of (maximal) almost disjoint families: very mad
families and (maximal) cofinitary groups. For the very mad families we prove
the basic existence results. We prove that MA implies there exist many pairwise
orthogonal families, and that CH implies that for any very mad family there is
one orthogonal to it. Finally we prove that the axiom of constructibility
implies that there exists a coanalytic very mad family. Cofinitary groups have
a natural action on the natural numbers. We prove that a maximal cofinitary
group cannot have infinitely many orbits under this action, but can have any
combination of any finite number of finite orbits and any finite (but nonzero)
number of infinite orbits.

 We also prove that there exists a maximal cofinitary group into which each
countable group embeds. This gives an example of a maximal cofinitary group
that is not a free group. We start the investigation into which groups have
cofinitary actions. The main result there is that it is consistent that
$\bigoplus_{\alpha \in \aleph_1} \mathbb{Z}_2$ has a cofinitary action.
 Concerning the complexity of maximal cofinitary groups we prove that they
cannot be $K_\sigma$, but that the axiom of constructibility implies that there
exists a coanalytic maximal cofinitary group. We prove that the least
cardinality $\mathfrak{a}_g$ of a maximal cofinitary group can consistently be
less than the cofinality of the symmetric group. Finally we prove that
$\mathfrak{a}_g$ can consistently be bigger than all cardinals in Cicho\'n's
diagram.}

\arXiv{0910.2318}
{Forcing, games and families of closed sets}
{Marcin Sabok}
{We propose a new, game-theoretic, approach to the idealized forcing, in terms
of fusion games. This generalizes the classical approach to the Sacks and the
Miller forcing. For definable ($\mathbf{\Pi}^1_1$ on $\mathbf{\Sigma}^1_1$)
$\sigma$-ideals we show that if a $\sigma$-ideal is generated by closed sets,
then it is generated by closed sets in all forcing extensions. We also prove an
infinite-dimensional version of the Solecki dichotomy for analytic sets. Among
examples, we investigate the $\sigma$-ideal $\E$ generated by closed null sets
and $\sigma$-ideals connected with not piecewise continuous functions.}

\arXiv{0910.3091}
{CH, a problem of Rolewicz and bidiscrete systems}
{Mirna Dzamonja, Istvan Juhasz}
{We give a construction under $CH$ of a non-metrizable compact Hausdorff space
$K$ such that any uncountable semi-biorthogonal sequence in $C(K)$ must be of a
very specific kind. The space $K$ has many nice properties, such as being
hereditarily separable, hereditarily Lindel\"of and a 2-to-1 continuous
preimage of a metric space, and all Radon measures on $K$ are separable.
However $K$ is not a Rosenthal compactum.

 We introduce the notion of bidiscrete systems in compact spaces and prove
that every infinite compact Hausdorff space $K$ must have a bidiscrete system
of size $d(K)$, the density of $K$. This, in particular, implies that $C(K)$
has a biorthogonal system of size $d(K)$ (known for $d(K)=\aleph_1$) and the
known result that $K^2$ has a discrete subspace of size $d(K)$.}

\arXiv{0910.4106}
{A Note on Monotonically Metacompact Spaces}
{Harold R. Bennett, Klaas Pieter Hart, David J. Lutzer}
{We show that any metacompact Moore space is
monotonically metacompact and use that result to characterize monotone
metacompactness in certain generalized ordered (GO)spaces.  We show, for
example, that a generalized ordered space with a
$\sigma$-closed-discrete dense subset is metrizable if and only if it
is monotonically (countably) metacompact, that a monotonically
(countably) metacompact GO-space is hereditarily paracompact, and that
a locally countably compact GO-space is metrizable if and only if it
is monotonically (countably) metacompact. We give an example of a
non-metrizable LOTS that is monotonically metacompact, thereby
answering a question posed by S. G.  Popvassilev. We also give
consistent examples showing that if there is a Souslin line, then
there is one Souslin line that is monotonically countable metacompact,
and another Souslin line that is not monotonically countably
metacompact.}

\arXiv{0910.4107}
{Covering dimension and finite-to-one maps}
{Klaas Pieter Hart, Jan van Mill}
{Hurewicz characterized the dimension of separable metrizable spaces
by means of finite-to-one maps.
We investigate whether this characterization also holds in the class
of compact $F$-spaces of weight~$\fc$.
Our main result is that, assuming the Continuum Hypothesis, an $n$-dimensional
compact $F$-space of weight~$\fc$ is the continuous image of a zero-dimensional
compact Hausdorff space by an at most $2^n$-to-$1$ map.}

\arXiv{0911.0145}{The Whyburn property in the class of $P$-spaces}
{Angelo Bella, Camillo Costantini and Santi Spadaro}
{We investigate the Whyburn and weakly Whyburn property in the class of
$P$-spaces, that is spaces where every $G_\delta$ set is open. We construct
examples of non-weakly Whyburn $P$-spaces of size continuum, thus giving a
negative answer under CH to a question of Pelant, Tkachenko, Tkachuk and
Wilson. In addition, we show that the weak Kurepa Hypothesis (a set-theoretic
assumption weaker than CH) implies the existence of a non-weakly Whyburn
$P$-space of size $\aleph_2$. Finally, we consider the behavior of the
above-mentioned properties under products; we show in particular that the
product of a Lindel\"of weakly Whyburn P-space and a Lindel\"of Whyburn
$P$-space is weakly Whyburn, and we give a consistent example of a non-Whyburn
product of two Lindel\"of Whyburn $P$-spaces.}

\arXiv{0911.0332}
{The Alaoglu theorem is equivalent to the Tychonoff theorem for compact Hausdorff spaces}
{Stefano Rossi}
{In this brief note we provide a simple approach to show that the Alaoglu
theorem and the Tychonoff theorem for compact Hausdorff spaces are equivalent.}

\arXiv{0911.2774}
{On splitting infinite-fold covers}
{M\'arton Elekes, Tam\'as M\'atrai, Lajos Soukup}
{Let $X$ be a set, $\kappa$ be a cardinal number and let $\cA$ be a family of
subsets of $X$ which covers each $x\in X$ at least $\kappa$ times. What
assumptions can ensure that $\cA$ can be decomposed into $\kappa$ many disjoint
subcovers?
 We examine this problem under various assumptions on the set $X$ and on the
cover $\cA$: among other situations, we consider covers of topological spaces
by closed sets, interval covers of linearly ordered sets and covers of
$\bbR^{n}$ by polyhedra and by arbitrary convex sets. We focus on these
problems mainly for infinite $\kappa$. Besides numerous positive and negative
results, many questions turn out to be independent of the usual axioms of set
theory.}

\arXiv{0911.3861}
{Maximality of ideal-independent sets}
{Corey Thomas Bruns}
{In this note we derive a property of maximal ideal-independent subsets of
boolean algebras which has corollaries regarding the continuum cardinals p and
$s_{mm}(P(\omega)/\mathrm{fin})$.}

\arXiv{0911.3833}
{On Galvin's lemma and Ramsey spaces}
{Jose G. Mijares}
{An abstract version of Galvin's lemma is proven, within the framework of the
theory of Ramsey spaces. Some instances of it are explored.}

\arXiv{0911.4075}
{The topological structure of (homogeneous) spaces and groups with
 countable cs$^*$-character}
{Taras Banakh and Lyubomyr Zdomskyy}
{In this paper we introduce and study three new cardinal topological
invariants called the cs$^*$-, cs-, and sb-characters. The class of topological
spaces with countable cs$^*$-character is closed under many topological operations
and contains all $\aleph$-spaces and all spaces with point-countable cs$^*$-network.
Our principal result states that each non-metrizable sequential topological
group with countable cs$^*$-character has countable pseudo-character and contains
an open $k_\omega$-subgroup.}

\arXiv{0911.4081}
{On topological groups containing a Fr\'echet-Urysohn fan}
{Taras Banakh}
{Suppose $G$ is a topological group containing a (closed) topological copy of
the Fr\'echet-Urysohn fan. If $G$ is a perfectly normal sequential space (a normal
$k$-space) then every closed metrizable subset in $G$ is locally compact.
Applying this result to topological groups whose underlying topological space
can be written as a direct limit of a sequence of closed metrizable subsets, we
get that every such a group either is metrizable or is homeomorphic to the
product of a $k_\omega$-space and a discrete space.}

\arXiv{0912.0497}
{Hewitt-Marczewski-Pondiczery type theorem for abelian groups and Markov's potential density}
{Dikran Dikranjan, Dmitri Shakhmatov}
{For an uncountable cardinal $\tau$ and a subset $S$ of an abelian group $G$, the following conditions are equivalent:
\begin{itemize}
\item[(i)] $|\{ns:s\in S\}|\ge \tau$ for all integers $n\ge 1$;
\item[(ii)] there exists a group homomorphism $\pi:G\to \mathbb{T}^{2^\tau}$ such that $\pi(S)$ is dense in  $\mathbb{T}^{2^\tau}$.
\end{itemize}
Moreover, if $|G|\le 2^{2^\tau}$, then the following item can be added to this list:
\begin{itemize}
\item[(iii)] there exists an isomorphism $\pi:G\to G'$ between $G$ and a subgroup $G'$ of $\mathbb{T}^{2^\tau}$ such that $\pi(S)$ is dense in  $\mathbb{T}^{2^\tau}$.
\end{itemize}
We prove that the following conditions are equivalent for an uncountable subset $S$ of an abelian group $G$ that is either (almost) torsion-free or divisible:
 \begin{itemize}
  \item[(a)] $S$ is $\mathcal{T}$-dense in $G$ for some Hausdorff group topology $\mathcal{T}$ on $G$;
  \item[(b)] $S$ is $\mathcal{T}$-dense in some precompact Hausdorff group topology $\mathcal{T}$ on $G$;
  \item[(c)] $|\{ns:s\in S\}|\ge \min\left\{\tau:|G|\le 2^{2^\tau}\right\}$ for every integer $n\ge 1$.
\end{itemize}
This partially resolves a question of Markov going back to 1946.}

\arXiv{0911.4978}
{Product between ultrafilters and applications to the Connes' embedding problem}
{V. Capraro, L. Paunescu}
{In this paper we want to apply the notion of product between ultrafilters to
answer several questions which arise around the Connes' embedding problem. For
instance, we will give a simplification and generalization of a theorem by
Radulescu; we will prove that ultraproduct of hyperlinear groups is still
hyperlinear and consequently the von Neumann algebra of the free group with
uncountable many generators is embeddable into $R^{\omega}$. This follows also
from a general construction that allows, starting from an hyperlinear group, to
find a family of hyperlinear groups. We will introduce the notion of
hyperlinear pair and we will use it to give some other characterizations of
hyperlinearity. We will prove also that the cross product of a hyperlinear
group via a profinite action is embeddable into $R^{\omega}$.}

\arXiv{0912.0406}
{A dichotomy for the number of ultrapowers}
{Ilijas Farah, Saharon Shelah}
{We prove a strong dichotomy for the number of ultrapowers of a given
countable model associated with nonprincipal ultrafilters on N. They are either
all isomorphic, or else there are $2^{2^{\aleph_0}}$ many nonisomorphic
ultrapowers. We prove the analogous result for metric structures, including
$C^*$-algebras and II$_1$ factors, as well as their relative commutants and
include several applications. We also show that the $C^*$-algebra $B(H)$ always has
nonisomorphic relative commutants in its ultrapowers associated with
nonprincipal ultrafilters on $\N$.}

\arXiv{0912.0431}
{Optimal Matrices of Partitions and an Application to Souslin Trees}
{Gido Scharfenberger-Fabian}
{The basic result of this note is a statement about the existence of families
of partitions of the set of natural numbers with some favourable properties,
the n-optimal matrices of partitions. We use this to improve a decomposition
result for strongly homogeneous Souslin trees. The latter is in turn applied to
separate strong notions of rigidity of Souslin trees, thereby answering a
considerable portion of a question of Fuchs and Hamkins.}

\arXiv{0912.2946}
{Locally minimal topological groups}
{Lydia Au\ss enhofer, Mar\'ia Jes\'us Chasco, Dikran Dikranjan, Xabier Dom\'inguez}
{A Hausdorff topological group $(G,\tau)$ is called locally minimal if there
exists a neighborhood $U$ of 0 in $\tau$ such that $U$ fails to be a
neighborhood of zero in any Hausdorff group topology on $G$ which is strictly
coarser than $\tau.$ Examples of locally minimal groups are all subgroups of
Banach-Lie groups, all locally compact groups and all minimal groups.
 Motivated by the fact that locally compact NSS groups are Lie groups, we
study the connection between local minimality and the NSS property,
establishing that under certain conditions, locally minimal NSS groups are
metrizable.
 \par
 A symmetric subset of an abelian group containing zero is said to be a GTG
set if it generates a group topology in an analogous way as convex and
symmetric subsets are unit balls for pseudonorms on a vector space. We consider
topological groups which have a neighborhood basis at zero consisting of GTG
sets. Examples of these locally GTG groups are: locally pseudo--convex spaces,
groups uniformly free from small subgroups (UFSS groups) and locally compact
abelian groups. The precise relation between these classes of groups is
obtained: a topological abelian group is UFSS if and only if it is locally
minimal, locally GTG and NSS. We develop a universal construction of GTG sets
in arbitrary non-discrete metric abelian groups, that generates a strictly
finer non-discrete UFSS topology and we characterize the metrizable abelian
groups admitting a strictly finer non-discrete UFSS group topology. We also
prove that a bounded abelian group $G$ admits a non-discrete locally minimal
and locally GTG group topology iff $|G|\geq \fc$.}

\ed